\documentclass[]{article}
\usepackage{lineno,epsfig,graphics,amssymb,amsmath,graphicx,pifont,soul,esint,float,multirow,array,epstopdf,color,url}

\newcommand{\bl}[1]{\boldsymbol{#1}}

\providecommand{\BibitemShut}[1]{}

\begin{document}

\title{A stabilized formulation for the solution of the incompressible unsteady Stokes equations in the frequency domain}

\author{Mahdi Esmaily\\Cornell University}
\date{} 
  
\maketitle

\begin{abstract}
   A stabilized finite element method is introduced for the simulation of time-periodic creeping flows, such as those found in the cardiorespiratory systems. The new technique, which is formulated in the frequency rather than time domain, strictly uses real arithmetics and permits the use of similar shape functions for pressure and velocity for ease of implementation. It involves the addition of the Laplacian of pressure to the continuity equation with a complex-valued stabilization parameter that is derived systematically from the momentum equation. The numerical experiments show the excellent accuracy and robustness of the proposed method in simulating flows in complex and canonical geometries for a wide range of conditions. The present method significantly outperforms a traditional solver in terms of both computational cost and scalability, which lowers the overall solution turnover time by several orders of magnitude.
\end{abstract}

\section{Introduction} 
\label{sec:intro}
Simulation of time-periodic creeping flows, such as cardiorespiratory flows in smaller vessels, can significantly benefit from a formulation that is expressed in the frequency rather than the time domain. 
Firstly, the boundary conditions in these problems, which typically vary smoothly in time, can be represented via a handful of Fourier modes. 
By solving for those few selected Fourier modes rather than integrating over thousands of time steps, a frequency formulation reduces the cost of a simulation by orders of magnitude. 
The absence of homogeneous solution in the frequency formulation presents a second major cost advantage. 
In a standard time formulation, the transient solution obtained at the beginning of the simulation has little significance, yet it must be computed before the particular solution can be obtained. 
Since a frequency formulation is independent of the initial conditions, it directly produces the particular solution, thus avoiding this costly and unnecessary computation. 
Thirdly, while the frequency formulation is embarrassingly parallelizable, a standard time formulation can be hardly parallelized in time. 
Thus, computations that are based on frequency formulation can be scaled to a much larger number of processors, permitting a much shorter solution turnover time. 
Fourthly, the time integration error present in the standard time formulation is absent in the frequency formulation. 
The truncation error associated with the selection of a finite number of Fourier modes is, however, present in the frequency formulation. 
Lastly, the stability consideration associated with the time integrator is no longer a concern in a frequency formulation as the solution is obtained from solving a boundary value problem.

The advantages of the frequency formulation for the solution of the incompressible unsteady Stokes equations enumerated above were shown in practice in a recent article \cite{meng2020time}. 
That included one to two orders of magnitude reduction in cost and improvement in scalability by the number of computed modes. 
Despite these attractive results, the formulation presented in \cite{meng2020time} had several shortcomings. 
Firstly, it relies on complex arithmetics and thus was hard to implement by requiring a significant change in the existing implementation of both the fluid and linear solver. 
Furthermore, a Bubnov-Galerkin formulation was employed in that case to satisfy the inf-sup condition  \cite{ladyzhenskaya1969mathematical, babuvska1971error, brezzi1974existence}. 
That led to two additional shortcomings, namely the requirement to use mixed shape functions for pressure and velocity and also the stiffness of the tangent matrix. 
The use of linear shape function for pressure and quadratic shape function for velocity is not a convenient choice, particularly in the case of complex geometries. 
The issue associated with the stiffness of the tangent matrix is caused by a zero block on the diagonal of the stiffness matrix (given that the continuity equation does not depend on the pressure), which delays the convergence of the iterative linear solver. 
It was shown that this slow convergence contributes to an order of magnitude increase in cost, particularly at higher modes where the condition number of the stiffness matrix increases.

The present article introduces a stabilized formulation to overcome the aforementioned issues. 
Namely, the new formulation uses purely real arithmetics, permits the use of equal order shape functions, and avoids the zero block in the tangent matrix for faster convergence of the linear solver. 
The article is organized as follows. 
In Section \ref{sec:method}, the stabilized formulation is presented. 
Then in Section \ref{sec:results}, the proposed formulation is tested using a canonical and a complex patient-specific geometry. 
Lastly, Section \ref{sec:conclusion} contains the concluding remarks. 

\section{A time-spectral stabilized formulation for Stokes equation}  \label{sec:method}
Creeping flows, such as those found in cardiorespiratory flows, can be modeled using the incompressible unsteady Stokes equation. 
Taking the temporal Fourier transformation of the unsteady Stokes equations produces a boundary value problem that is stated as 
\begin{equation}
   \begin{alignedat}{2}
      \hat{j} \rho\omega \bl u &= - \nabla p + \nabla \cdot (\mu \nabla \bl u)  \; \; \; && \mathrm {in} \;\Omega, \\
        \nabla \cdot \bl u &= 0 && \mathrm {in} \; \Omega, \\
        \bl u &= \bl g && \mathrm {on} \; \Gamma_{\rm g},\\
        (- p \bl I + \mu\nabla \bl u)\cdot \bl n &= \bl h && \mathrm {on} \; \Gamma_{\rm h}.
    \end{alignedat}
    \label{stokes}
\end{equation}
where $\hat j = \sqrt{-1}$, $\omega$ is the oscillation frequency, $\bl x$ is position, $\bl u(\omega, \bl x)$ is the velocity, $p(\omega,\bl x)$ is pressure, $\bl g(\omega,\bl x)$ is the imposed velocity on the Dirichlet boundary $\Gamma_{\rm g}$, and $\bl h(\omega,\bl x)$ is the imposed traction on the Neumann boundary $\Gamma_{\rm h}$. 
Note that the time derivative term in the Stokes equations is transformed to a complex-valued source term in Eq. \eqref{stokes}.
Thus, the unsteady Stokes equations, when expressed in the frequency domain, behave similarly to the steady Stokes equations with a nonzero complex source term. 

Since $\bl u$ and $\bl p$ are complex-valued variables in general, numerical simulation of Eq. \eqref{stokes} will require complex arithmetics \cite{meng2020time}. 
For a real-valued formulation, the state variables as well as boundary conditions can be expressed as 
\begin{equation}
   \begin{split}
   \bl u &= \bl u_r + \hat j \bl u_i,\\
   p &= p_r + \hat j p_i,\\
   \bl g &= \bl g_r + \hat j \bl g_i,\\
   \bl h &= \bl h_r + \hat j \bl h_i. 
   \end{split}
   \label{decomp}
\end{equation}
Using these changes of variables, Eq. \eqref{stokes} can be rewritten in the real domain as 
\begin{equation}
   \begin{alignedat}{3}
      - \rho\omega \bl u_i &= - \nabla p_r + \nabla \cdot (\mu \nabla \bl u_r), & \rho \omega \bl u_r &= - \nabla p_i + \nabla \cdot (\mu \nabla \bl u_i) \;\;\;&& \mathrm {in} \;\Omega, \\
        \nabla \cdot \bl u_r &= 0, & \nabla \cdot \bl u_i &= 0 && \mathrm {in} \; \Omega, \\
        \bl u_r &= \bl g_r, & \bl u_i &= \bl g_i  && \mathrm {on} \; \Gamma_{\rm g},\\
        (- p_r \bl I + \mu\nabla \bl u_r)\cdot \bl n &= \bl h_r, & (- p_i \bl I + \mu\nabla \bl u_i)\cdot \bl n &= \bl h_i && \mathrm {on} \; \Gamma_{\rm h}.\\
\end{alignedat}
    \label{stokes_split}
\end{equation}

To weak form of Eq. \eqref{stokes_split} can be expressed as finding $\bl u_r \in \mathcal {\bl S}_r$, $\bl u_i \in \mathcal {\bl S}_i$ and $p_r, p_i \in L^2$ such that for any $\bl w_r, \bl w_i \in \bl{\mathcal W}$ and $q_r, q_i \in L^2$ the following must hold
\begin{equation}
    \begin{split}
       B_{\rm G} =  \int_\Omega\bigg[ &-\rho \omega \bl w_r \cdot \bl u_i 
                                       + \nabla \bl w_r : ( -p_r I + \mu \nabla \bl u_r) 
                                       - q_r \nabla \cdot \bl u_r\\
                                      &-\rho\omega \bl w_i \cdot \bl u_r 
                                       - \nabla \bl w_i : ( -p_i I + \mu \nabla \bl u_i) 
                                      + q_i \nabla \cdot \bl u_i \bigg] \rm d \Omega, \\
        F_{\rm G} = \int_{\Gamma_{\rm h}}\bigg[ &\bl w_r \cdot \bl h_r - \bl w_i \cdot \bl h_i \bigg] \rm d \Gamma.
    \end{split}
    \label{weak_g}
\end{equation}
In this equation, $\bl w_r$, $\bl w_i$, $q_r$ and $q_i$ are test functions for velocity and pressure and 
\begin{equation}
    \begin{split}
        \bl{\mathcal S}_r  &= \left\{ \bl u_r | \bl u_r \in(H^1)^{n_{\rm sd}} ,\; {\bl u = \bl g_r} \;\rm{on}\; \Gamma_{\rm g} \right\}, \\
        \bl{\mathcal S}_i  &= \left\{ \bl u_i | \bl u_i \in(H^1)^{n_{\rm sd}} ,\; {\bl u = \bl g_i} \;\rm{on}\; \Gamma_{\rm g} \right\}, \\
        \bl{\mathcal W} &= \left\{ \bl w | \bl w \in(H^1)^{n_{\rm sd}} ,\; \bl w = \bl 0 \; \rm{on}\; \Gamma_{\rm g} \right\}.
    \end{split}
\end{equation}
In above, $L^2$ denotes the space of scalar-valued functions that are square-integrable on $\Omega$.
Also, $(H^1)^{n_{sd}}$ denotes the space of vector-valued functions with square-integrable derivatives on $\Omega$.

In obtaining $B_{\rm G}$ and $F_{\rm G}$ in Eq. \eqref{weak_g}, the imaginary part of the momentum equation and real part of the continuity equation were multiplied by $-1$. 
The sign of those equations was changed to ensure the tangent matrix remains symmetric. 

To stabilize Eq. \eqref{weak_g} for equal order shape functions for $\bl u$ and $p$ and relate the continuity equation to pressure, the Laplacian of pressure is added to the continuity equation. 
To systematically derive this term, consider the divergence of the momentum equation from Eq. \eqref{stokes} that is
\begin{equation}
   \hat{j}\rho \omega  \nabla \cdot \bl u = - \nabla^2 p + \nabla \cdot \left[\nabla \cdot (\mu \nabla \bl u) \right]. 
    \label{stokes_S1}
\end{equation}
Since the last term in Eq. \eqref{stokes_S1} involves the third derivative of velocity, it will vanish in the interior of an element when linear shape functions are employed. 
It is thus approximated using the characteristic size of the element $H$ as $\nabla \cdot \left[\nabla \cdot (\mu \nabla \bl u) \right] \approx (\mu/H^2) \nabla \cdot \bl u $. 
Thus, Eq. \eqref{stokes_S1} can be written as
\begin{equation}
   \nabla \cdot \bl u - \tau \nabla^2 p = 0,
    \label{stokes_S2}
\end{equation}
where
\begin{equation}
   \tau = \frac{1}{(\mu/H^2) - \rho\omega \hat j}.
    \label{tau_derived}
\end{equation}
The numerical experiments show that incorporating a small constant in the definition of $\tau$ improves the solution accuracy and reduces the number of linear solver iterations. 
They also show $\bl \xi \in \mathbb R^{n_{\rm sd}} \times  \mathbb R^{n_{\rm sd}}$, which is the covariant tensor obtained from a mapping between the physical and parent elements, provides a good approximation for $H$ \cite{bazilevs2013computational}. 
Thus, in practice, the real and imaginary component of $\tau$ are defined as 
\begin{equation}
   \tau_r = \frac{c\mu \sqrt{\bl \xi:\bl \xi}}{(\rho\omega)^2+ \mu^2\bl \xi :\bl \xi } ,\;\;\;\; \tau_i = \frac{c\rho\omega }{(\rho\omega)^2 + \mu^2 \bl \xi :\bl \xi}.
    \label{tau}
\end{equation}
The numerical simulation involving 2D triangular elements and 3D tetrahedral elements show $c=2^{-5} \approx 0.03$ produces satisfactory results. 
This value of $c$ is utilized for all the computations reported in Section \ref{sec:results}. 

From Eqs. \eqref{stokes_S2} and \eqref{tau}, the stabilization terms are computed as
\begin{equation}
   B_{\rm S} = \sum_e \int_{\Omega_e} \bigg[ -\nabla q_r \cdot \left(\tau_r \nabla p_r - \tau_i \nabla p_i\right) +  \nabla q_i \cdot \left(\tau_r \nabla p_i + \tau_i \nabla p_r \right) \bigg] \rm d \Omega,
    \label{weak_s}
\end{equation}
where the signs are selected to be consistent with those in Eq. \eqref{weak_g}. 
These integrals are added to the Galerkin's weak form from Eq. \eqref{weak_g} to obtain 
\begin{equation}
   B_{\rm G} + B_{\rm S} = F_{\rm G}. 
   \label{weak}
\end{equation}

Derivation of the matrix form of Eq. \eqref{weak} follows a standard process and, thus, not included here in detail. 
The result is 
\begin{equation}
   \left[ \begin{array}{cccc}
      \mu L_{AB} \bl \delta & -\bl G_{AB} & -\rho \omega M_{AB} \bl \delta & \bl 0 \\ [0.5cm]
      -\bl D_{AB}  & -\tau_r L_{AB} & \bl 0 & \tau_i L_{AB} \\  [0.5cm]
      -\rho \omega M_{AB} \bl \delta  & \bl 0  & -\mu L_{AB}\bl \delta & \bl G_{AB}   \\  [0.5cm]
      \bl 0 & \tau_i L_{AB} & \bl D_{AB}  & \tau_r L_{AB}  \\ 
   \end{array}
   \right] \left[  \begin{array}{c}
      \bl U_{rB} \\ [0.5cm] P_{rB} \\ [0.5cm] \bl U_{iB} \\ [0.5cm] P_{iB}
   \end{array}\right] = -\left[  \begin{array}{c} 
      \bl R^m_{rA} \\ [0.5cm] R^c_{rA} \\ [0.5cm] \bl R^m_{iA} \\ [0.5cm] R^c_{iA}
   \end{array}\right]
   \label{axb}
\end{equation}
where 
\begin{equation}
   \begin{split}
      L_{AB} &= \int_\Omega \nabla N_A\cdot \nabla N_B  {\rm d} \Omega, \\
      \bl G_{AB} &= \int_\Omega \nabla N_A N_B {\rm d} \Omega, \\
      \bl D_{AB} &= \int_\Omega N_A \nabla N_B {\rm d} \Omega, \\
      M_{AB} &= \int_\Omega N_A N_B {\rm d} \Omega, \\
      \end{split}
   \label{ab_def}
\end{equation}
and 
\begin{equation}
   \begin{split}
      \bl R^m_{rA} &= -\int_{\Gamma_{\rm h}} N_A\bl h_r {\rm d}\Gamma + \mu L_{AB} \bl g_{rB} - \rho \omega M_{AB} \bl g_{iB}, \\
      \bl R^m_{iA} &= \int_{\Gamma_{\rm h}} N_A\bl h_i {\rm d}\Gamma - \mu L_{AB} \bl g_{iB} - \rho \omega M_{AB} \bl g_{rB}, \\
       R^c_{rA} &= -\bl D_{AB} \bl g_{rB}, \\ 
       R^c_{iA} &= \bl D_{AB} \bl g_{iB}. 
      \end{split}
   \label{r_def}
\end{equation}
In Eq. \eqref{axb}, $\bl U_{rB}$ and $\bl U_{iB}$ are the real and imaginary component of velocity vector at node $B$, respectively, and $P_{rB}$ and $P_{iB}$ are the real and imaginary component of pressure at node $B$, respectively. 

\emph{Remarks on Eq. \eqref{axb}} 
\begin{enumerate}
\item The tangent matrix in Eq. \eqref{axb} is independent of the solution, thus, there is no need for Newton-Raphson iterations when this scheme is implemented. 
   This linear property, which is a result of the linearity of the Stokes equations and that of the designed stabilization scheme, permits one to obtain the final solution via a single linear solution. 
\item There is a linear relationship between the number of nonzero blocks in the tangent matrix and the number of spatial dimensions $n_{\rm sd}$. 
   In total, the tangent matrix contains $4(n_{\rm sd} + 1)^2$ blocks from which only $8n_{\rm sd} + 4$ are nonzero.
That translates to less than half of all blocks for 3D problems (28 out of 64 are nonzero). 
Despite its relatively large size, this matrix is symmetric and can be solved using an efficient iterative solver. 
\item It is possible to design a linear solver with a specialized library for the matrix-vector product that only operates on the nonzero blocks of the tangent matrix. 
   Such an implementation is expected to roughly reduce the cost of 3D computations by half. 
      What is presented below, however, does not take advantage of this optimization and is based on a standard matrix-vector product library that is developed in-house \cite{Esmaily2015DS}. 
\item A set of iterative techniques are tested for solving this linear system including, conjugate gradient, successive over-relaxation, bi-conjugate gradient, generalized minimal residual \cite{saad1986gmres}, and bi-partitioned methods \cite{Esmaily2015BIPN}.
      A symmetric Jacobi preconditioner was used for all these cases, producing diagonal entries with a mod of one. 
      The results show that the conjugate gradient is the most efficient technique for solving the linear system in Eq. \eqref{axb} and thus used for all the cases presented in Section \ref{sec:results}. 
      This superior performance is despite the nonmonotonic convergence of the conjugate gradient and the fact that the tangent matrix is symmetric but indefinite. 
      Nevertheless, the conjugate gradient successfully converges for all cases considered below while being the least costly method among all techniques enumerated above.  
\item The introduced scheme only requires the use of a single complex-valued stabilization parameter, i.e. $\tau$ defined in Eq. \eqref{tau}. 
   There is an arbitrary constant $c$ that is incorporated into the definition of $\tau$. 
    The numerical results involving 2D and 3D elements show that using a different value for $c$ has a minimal effect on the accuracy of this stabilized formulation.
      Also, the numerical experiments show this formulation is fairly robust if the element length scale in Eq. \eqref{tau} is based on a parameter other than $\bl \xi$ (e.g., the Jacobian of element mapping or its volume).
\item The imaginary component of $\tau$, namely $\tau_i$, is crucially necessary at high frequencies. 
   This scheme will struggle to converge for $\omega H^2/\nu \gg 1$ if one ignores the contribution of the off-diagonal terms associated with $\tau_i$ in the tangent matrix. 
      Similarly at the steady state limit where $\omega \to 0$, $\tau_i \to 0$ and $\tau_r \to c/(\mu \sqrt{\bl \xi:\bl \xi})$. 
      At this limit, the form of $\tau_r$ becomes identical to the existing stabilization parameters that is utilized in the residual-based variational multiscale method at the limit of $\Delta t \to \infty$ and $\|\bl u\| \to 0$ \cite{shakib1991new,bazilevs2007variational}. 
   \item The real and imaginary unknowns are coupled via two sets of off-diagonal blocks in the tangent matrix, both of which are proportional to $\omega$. 
      As $\omega \to 0$, the two sets of unknowns become decoupled. 
      That is when $\bl u_r$ and $p_r$ depend only on $\bl h_r$ and $\bl g_r$ and not on $\bl h_i$ and $\bl g_i$ and vice versa. 
      This limit is physically known as the quasi steady limit where the acceleration term in the Stokes equations is negligible compared to the viscous and pressure terms. 
   \item Once the velocity and pressure unknowns are computed in the frequency domain, their temporal counterparts can be simply computed using Eq. \eqref{decomp} and
      \begin{equation}
         \begin{split}
            \hat {\bl u}(\bl x,t) &= \sum_{\omega} \bl u(\bl x,\omega) e^{\hat j\omega t},\\
            \hat p(\bl x,t) &= \sum_{\omega} p(\bl x,\omega) e^{\hat j\omega t}.
         \end{split}
      \end{equation}
\end{enumerate}

\section{Results} \label{sec:results}
Two sets of tests cases are considered in this section: 
1) a canonical case of an unsteady pipe flow, where the analytical Womersley solution is available for establishing the accuracy of the proposed scheme, and 
2) a patient-specific geometry for evaluating the performance of this method on more complex geometries. 
Although not presented here for the sake of brevity, simulation involving 2D channel flow has also been performed, where the results are in line with the above 3D cases. 

The method described above is implemented in an in-house finite element solver.
This solver is parallelized using a message passing interface (MPI).
The workload is only parallelized using spatial partitioning by employing ParMETIS library \cite{METIS}. 
Further parallelization across different frequencies is not considered here as it involves a trivial process of running a series of simulations with different $\omega$. 
All computations are performed on a cluster of AMD Opteron$^{\rm TM}$ 6378 processors that are interconnected via a QDR Infiniband.

Unless stated otherwise, a tolerance of $\epsilon_{LS} = 10^{-3}$ is used for the conjugate gradient to solve the linear systems. 

\subsection{Oscillatory pipe flow} 
An oscillatory laminar pipe flow is considered for the first test case. 
A pipe with a length to radius ratio of $L/R=15$ is considered with an oscillatory unit inlet and zero outlet Neumann boundary condition (i.e., $\bl h_r = 1 \bl n$ and $\bl h_i =\bl 0$ on the inlet and $\bl h=\bl 0$ on the outlet). 
The oscillation frequency $\omega$ is varied to simulate flow at eleven Womersley numbers $\alpha=R \sqrt{\rho \omega/\mu} = 0, \sqrt 2, 2, \cdots, 2^5$. 
Three tetrahedral meshes (M1, M2, and M3) are utilized for spatial discretization (Table \ref{table:msh}).
All these computations are performed using 16 processors unless stated otherwise.

\begin{table}[H]
  \centering
   \caption{Tetrahedral meshes used for discretization of the 3D pipe flow (M1--M3). 
  $N_{\rm ele}$ and $N_{\rm nds}$ denote the numbers of elements and nodes, respectively. }
  \label{table:msh}
    \begin{tabular}{c|ccc} 
       \hline\hline
       & M1 & M2 & M3  \\\hline
       $N_{\rm ele}$ & 24,450 & 207,063 & 728,922  \\
       $N_{\rm nds}$ & 5,462 & 37,401 & 122,291  \\ \hline\hline
    \end{tabular}
\end{table}

An analytical solution is available for an oscillatory flow in a pipe, that is expressed in the frequency domain as \cite{Womersley1955method} 
\begin{equation}
   u_{\rm ref}(r,\alpha) = \left\{ \begin{array}{lr} 
    \displaystyle \frac{h}{4\mu L}(R^2-r^2), &  \alpha = 0, \\ 
    & \\
       \displaystyle -\frac{\hat{j}h R^2}{L\mu\alpha^2} \left[1- J_0(\hat j^{\frac3 2}\alpha)^{-1} J_0(\hat j^{\frac3 2} \alpha \frac{r}{R})\right],  & \alpha \neq 0,
    \end{array}\right.
    \label{wom_sol}
\end{equation}
where $J_0$ is the zero order Bessel function of the first kind and $h$ is the magnitude of the imposed Neumann boundary condition, which is one in this case. 
The Womersley solution from Eq. \eqref{wom_sol} is used as the reference solution to evaluate the performance of the proposed solver. 

\begin{figure}
  \centering
   \begingroup
  \makeatletter
  \providecommand\color[2][]{%
    \GenericError{(gnuplot) \space\space\space\@spaces}{%
      Package color not loaded in conjunction with
      terminal option `colourtext'%
    }{See the gnuplot documentation for explanation.%
    }{Either use 'blacktext' in gnuplot or load the package
      color.sty in LaTeX.}%
    \renewcommand\color[2][]{}%
  }%
  \providecommand\includegraphics[2][]{%
    \GenericError{(gnuplot) \space\space\space\@spaces}{%
      Package graphicx or graphics not loaded%
    }{See the gnuplot documentation for explanation.%
    }{The gnuplot epslatex terminal needs graphicx.sty or graphics.sty.}%
    \renewcommand\includegraphics[2][]{}%
  }%
  \providecommand\rotatebox[2]{#2}%
  \@ifundefined{ifGPcolor}{%
    \newif\ifGPcolor
    \GPcolortrue
  }{}%
  \@ifundefined{ifGPblacktext}{%
    \newif\ifGPblacktext
    \GPblacktexttrue
  }{}%
  \let\gplgaddtomacro\g@addto@macro
  \gdef\gplbacktext{}%
  \gdef\gplfronttext{}%
  \makeatother
  \ifGPblacktext
    \def\colorrgb#1{}%
    \def\colorgray#1{}%
  \else
    \ifGPcolor
      \def\colorrgb#1{\color[rgb]{#1}}%
      \def\colorgray#1{\color[gray]{#1}}%
      \expandafter\def\csname LTw\endcsname{\color{white}}%
      \expandafter\def\csname LTb\endcsname{\color{black}}%
      \expandafter\def\csname LTa\endcsname{\color{black}}%
      \expandafter\def\csname LT0\endcsname{\color[rgb]{1,0,0}}%
      \expandafter\def\csname LT1\endcsname{\color[rgb]{0,1,0}}%
      \expandafter\def\csname LT2\endcsname{\color[rgb]{0,0,1}}%
      \expandafter\def\csname LT3\endcsname{\color[rgb]{1,0,1}}%
      \expandafter\def\csname LT4\endcsname{\color[rgb]{0,1,1}}%
      \expandafter\def\csname LT5\endcsname{\color[rgb]{1,1,0}}%
      \expandafter\def\csname LT6\endcsname{\color[rgb]{0,0,0}}%
      \expandafter\def\csname LT7\endcsname{\color[rgb]{1,0.3,0}}%
      \expandafter\def\csname LT8\endcsname{\color[rgb]{0.5,0.5,0.5}}%
    \else
      \def\colorrgb#1{\color{black}}%
      \def\colorgray#1{\color[gray]{#1}}%
      \expandafter\def\csname LTw\endcsname{\color{white}}%
      \expandafter\def\csname LTb\endcsname{\color{black}}%
      \expandafter\def\csname LTa\endcsname{\color{black}}%
      \expandafter\def\csname LT0\endcsname{\color{black}}%
      \expandafter\def\csname LT1\endcsname{\color{black}}%
      \expandafter\def\csname LT2\endcsname{\color{black}}%
      \expandafter\def\csname LT3\endcsname{\color{black}}%
      \expandafter\def\csname LT4\endcsname{\color{black}}%
      \expandafter\def\csname LT5\endcsname{\color{black}}%
      \expandafter\def\csname LT6\endcsname{\color{black}}%
      \expandafter\def\csname LT7\endcsname{\color{black}}%
      \expandafter\def\csname LT8\endcsname{\color{black}}%
    \fi
  \fi
    \setlength{\unitlength}{0.0500bp}%
    \ifx\gptboxheight\undefined%
      \newlength{\gptboxheight}%
      \newlength{\gptboxwidth}%
      \newsavebox{\gptboxtext}%
    \fi%
    \setlength{\fboxrule}{0.5pt}%
    \setlength{\fboxsep}{1pt}%
\begin{picture}(9360.00,9360.00)%
    \gplgaddtomacro\gplbacktext{%
      \csname LTb\endcsname
      \put(991,6739){\makebox(0,0)[r]{\strut{}$-0.2$}}%
      \put(991,7176){\makebox(0,0)[r]{\strut{}$0$}}%
      \put(991,7612){\makebox(0,0)[r]{\strut{}$0.2$}}%
      \put(991,8049){\makebox(0,0)[r]{\strut{}$0.4$}}%
      \put(991,8486){\makebox(0,0)[r]{\strut{}$0.6$}}%
      \put(991,8922){\makebox(0,0)[r]{\strut{}$0.8$}}%
      \put(991,9359){\makebox(0,0)[r]{\strut{}$1$}}%
      \put(1123,6519){\makebox(0,0){\strut{}$0$}}%
      \put(1872,6519){\makebox(0,0){\strut{}$0.2$}}%
      \put(2620,6519){\makebox(0,0){\strut{}$0.4$}}%
      \put(3369,6519){\makebox(0,0){\strut{}$0.6$}}%
      \put(4117,6519){\makebox(0,0){\strut{}$0.8$}}%
      \put(4866,6519){\makebox(0,0){\strut{}$1$}}%
      \put(4492,9097){\makebox(0,0)[l]{\strut{}(a)}}%
      \put(2620,8573){\makebox(0,0)[l]{\strut{}$u_r^*$}}%
      \put(2620,7525){\makebox(0,0)[l]{\strut{}$u_i^*$}}%
    }%
    \gplgaddtomacro\gplfronttext{%
      \csname LTb\endcsname
      \put(452,8049){\rotatebox{-270}{\makebox(0,0){\strut{}$u_r^*;\; u_i^*$}}}%
    }%
    \gplgaddtomacro\gplbacktext{%
      \csname LTb\endcsname
      \put(5484,6739){\makebox(0,0)[r]{\strut{}$-0.6$}}%
      \put(5484,7113){\makebox(0,0)[r]{\strut{}$-0.4$}}%
      \put(5484,7488){\makebox(0,0)[r]{\strut{}$-0.2$}}%
      \put(5484,7862){\makebox(0,0)[r]{\strut{}$0$}}%
      \put(5484,8236){\makebox(0,0)[r]{\strut{}$0.2$}}%
      \put(5484,8610){\makebox(0,0)[r]{\strut{}$0.4$}}%
      \put(5484,8985){\makebox(0,0)[r]{\strut{}$0.6$}}%
      \put(5484,9359){\makebox(0,0)[r]{\strut{}$0.8$}}%
      \put(5616,6519){\makebox(0,0){\strut{}$0$}}%
      \put(6365,6519){\makebox(0,0){\strut{}$0.2$}}%
      \put(7113,6519){\makebox(0,0){\strut{}$0.4$}}%
      \put(7862,6519){\makebox(0,0){\strut{}$0.6$}}%
      \put(8610,6519){\makebox(0,0){\strut{}$0.8$}}%
      \put(9359,6519){\makebox(0,0){\strut{}$1$}}%
      \put(8985,9097){\makebox(0,0)[l]{\strut{}(b)}}%
    }%
    \gplgaddtomacro\gplfronttext{%
    }%
    \gplgaddtomacro\gplbacktext{%
      \csname LTb\endcsname
      \put(991,3556){\makebox(0,0)[r]{\strut{}$-0.35$}}%
      \put(991,3858){\makebox(0,0)[r]{\strut{}$-0.3$}}%
      \put(991,4159){\makebox(0,0)[r]{\strut{}$-0.25$}}%
      \put(991,4461){\makebox(0,0)[r]{\strut{}$-0.2$}}%
      \put(991,4762){\makebox(0,0)[r]{\strut{}$-0.15$}}%
      \put(991,5064){\makebox(0,0)[r]{\strut{}$-0.1$}}%
      \put(991,5365){\makebox(0,0)[r]{\strut{}$-0.05$}}%
      \put(991,5667){\makebox(0,0)[r]{\strut{}$0$}}%
      \put(991,5968){\makebox(0,0)[r]{\strut{}$0.05$}}%
      \put(991,6270){\makebox(0,0)[r]{\strut{}$0.1$}}%
      \put(1123,3336){\makebox(0,0){\strut{}$0$}}%
      \put(1872,3336){\makebox(0,0){\strut{}$0.2$}}%
      \put(2620,3336){\makebox(0,0){\strut{}$0.4$}}%
      \put(3369,3336){\makebox(0,0){\strut{}$0.6$}}%
      \put(4117,3336){\makebox(0,0){\strut{}$0.8$}}%
      \put(4866,3336){\makebox(0,0){\strut{}$1$}}%
      \put(4492,3827){\makebox(0,0)[l]{\strut{}(c)}}%
    }%
    \gplgaddtomacro\gplfronttext{%
      \csname LTb\endcsname
      \put(320,4913){\rotatebox{-270}{\makebox(0,0){\strut{}$u_r^*;\; u_i^*$}}}%
      \put(1783,5243){\makebox(0,0)[r]{\strut{}Ref.}}%
      \csname LTb\endcsname
      \put(1783,5023){\makebox(0,0)[r]{\strut{}M1}}%
      \csname LTb\endcsname
      \put(1783,4803){\makebox(0,0)[r]{\strut{}M2}}%
      \csname LTb\endcsname
      \put(1783,4583){\makebox(0,0)[r]{\strut{}M3}}%
    }%
    \gplgaddtomacro\gplbacktext{%
      \csname LTb\endcsname
      \put(5484,3556){\makebox(0,0)[r]{\strut{}$-0.07$}}%
      \put(5484,3827){\makebox(0,0)[r]{\strut{}$-0.06$}}%
      \put(5484,4099){\makebox(0,0)[r]{\strut{}$-0.05$}}%
      \put(5484,4370){\makebox(0,0)[r]{\strut{}$-0.04$}}%
      \put(5484,4642){\makebox(0,0)[r]{\strut{}$-0.03$}}%
      \put(5484,4913){\makebox(0,0)[r]{\strut{}$-0.02$}}%
      \put(5484,5184){\makebox(0,0)[r]{\strut{}$-0.01$}}%
      \put(5484,5456){\makebox(0,0)[r]{\strut{}$0$}}%
      \put(5484,5727){\makebox(0,0)[r]{\strut{}$0.01$}}%
      \put(5484,5999){\makebox(0,0)[r]{\strut{}$0.02$}}%
      \put(5484,6270){\makebox(0,0)[r]{\strut{}$0.03$}}%
      \put(5616,3336){\makebox(0,0){\strut{}$0$}}%
      \put(6365,3336){\makebox(0,0){\strut{}$0.2$}}%
      \put(7113,3336){\makebox(0,0){\strut{}$0.4$}}%
      \put(7862,3336){\makebox(0,0){\strut{}$0.6$}}%
      \put(8610,3336){\makebox(0,0){\strut{}$0.8$}}%
      \put(9359,3336){\makebox(0,0){\strut{}$1$}}%
      \put(8985,3827){\makebox(0,0)[l]{\strut{}(d)}}%
    }%
    \gplgaddtomacro\gplfronttext{%
    }%
    \gplgaddtomacro\gplbacktext{%
      \csname LTb\endcsname
      \put(991,468){\makebox(0,0)[r]{\strut{}$-0.02$}}%
      \put(991,905){\makebox(0,0)[r]{\strut{}$-0.015$}}%
      \put(991,1341){\makebox(0,0)[r]{\strut{}$-0.01$}}%
      \put(991,1778){\makebox(0,0)[r]{\strut{}$-0.005$}}%
      \put(991,2215){\makebox(0,0)[r]{\strut{}$0$}}%
      \put(991,2651){\makebox(0,0)[r]{\strut{}$0.005$}}%
      \put(991,3088){\makebox(0,0)[r]{\strut{}$0.01$}}%
      \put(1123,248){\makebox(0,0){\strut{}$0$}}%
      \put(1872,248){\makebox(0,0){\strut{}$0.2$}}%
      \put(2620,248){\makebox(0,0){\strut{}$0.4$}}%
      \put(3369,248){\makebox(0,0){\strut{}$0.6$}}%
      \put(4117,248){\makebox(0,0){\strut{}$0.8$}}%
      \put(4866,248){\makebox(0,0){\strut{}$1$}}%
      \put(1310,2826){\makebox(0,0)[l]{\strut{}(e)}}%
    }%
    \gplgaddtomacro\gplfronttext{%
      \csname LTb\endcsname
      \put(254,1778){\rotatebox{-270}{\makebox(0,0){\strut{}$u_r^*;\; u_i^*$}}}%
      \put(2994,-82){\makebox(0,0){\strut{}$r/R$}}%
    }%
    \gplgaddtomacro\gplbacktext{%
      \csname LTb\endcsname
      \put(5484,468){\makebox(0,0)[r]{\strut{}$-0.005$}}%
      \put(5484,842){\makebox(0,0)[r]{\strut{}$-0.004$}}%
      \put(5484,1217){\makebox(0,0)[r]{\strut{}$-0.003$}}%
      \put(5484,1591){\makebox(0,0)[r]{\strut{}$-0.002$}}%
      \put(5484,1965){\makebox(0,0)[r]{\strut{}$-0.001$}}%
      \put(5484,2339){\makebox(0,0)[r]{\strut{}$0$}}%
      \put(5484,2714){\makebox(0,0)[r]{\strut{}$0.001$}}%
      \put(5484,3088){\makebox(0,0)[r]{\strut{}$0.002$}}%
      \put(5616,248){\makebox(0,0){\strut{}$0$}}%
      \put(6365,248){\makebox(0,0){\strut{}$0.2$}}%
      \put(7113,248){\makebox(0,0){\strut{}$0.4$}}%
      \put(7862,248){\makebox(0,0){\strut{}$0.6$}}%
      \put(8610,248){\makebox(0,0){\strut{}$0.8$}}%
      \put(9359,248){\makebox(0,0){\strut{}$1$}}%
      \put(5803,2826){\makebox(0,0)[l]{\strut{}(f)}}%
    }%
    \gplgaddtomacro\gplfronttext{%
      \csname LTb\endcsname
      \put(7487,-82){\makebox(0,0){\strut{}$r/R$}}%
    }%
    \gplbacktext
    \put(0,0){\includegraphics[width={468.00bp},height={468.00bp}]{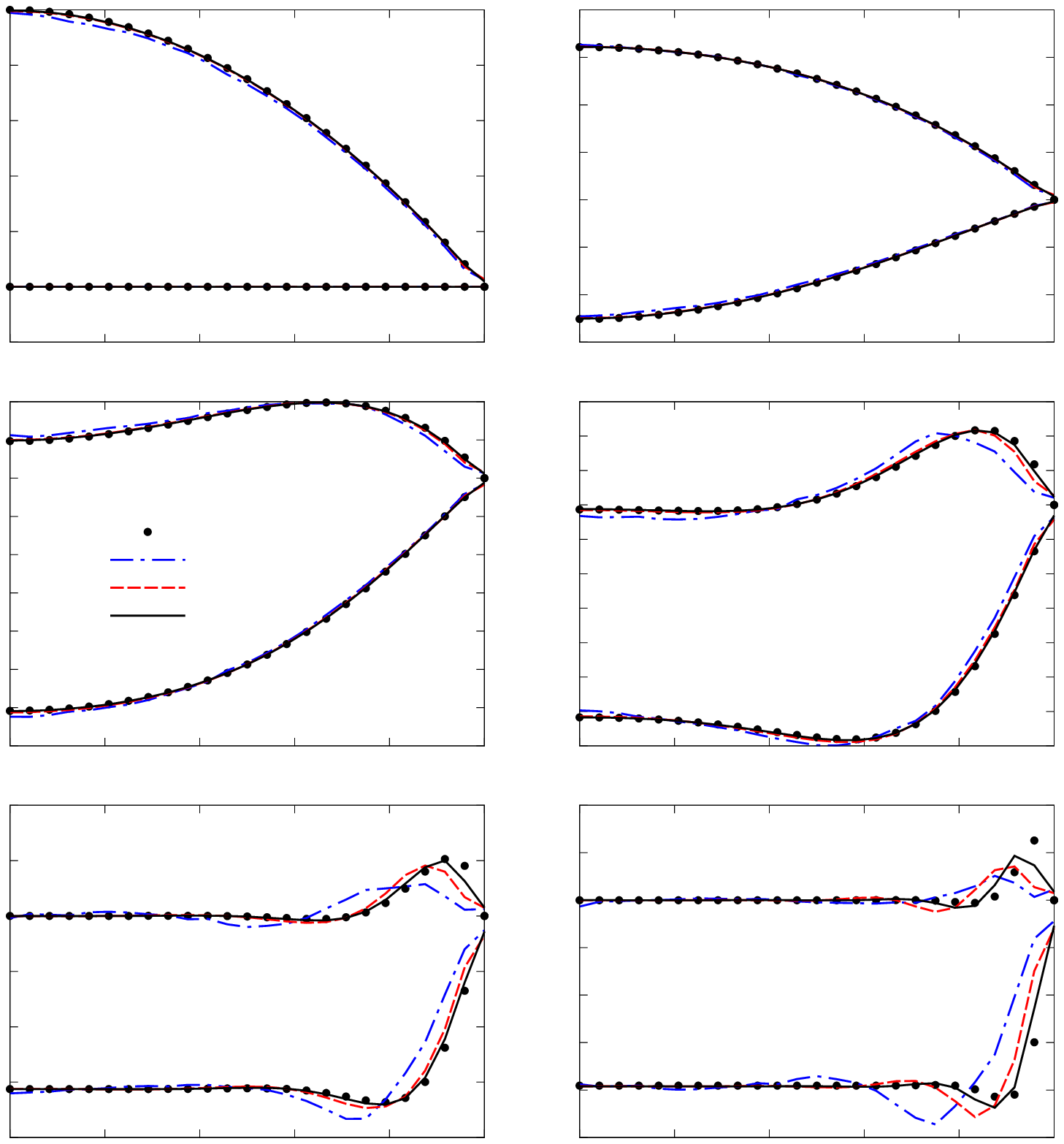}}%
    \gplfronttext
  \end{picture}%
\endgroup
   \caption{Axial velocity of a pulsating pipe flow as a function of radius computed from the proposed formulation on three meshes (Table \ref{table:msh} where M1, M2, and M3 are represented by dash-dot, dashed, and solid lines, respectively) and the analytical reference solution (Eq. \eqref{wom_sol} that is shown by dots).
   Different panels correspond to different Womersley number $\alpha= R\sqrt{\rho\omega/\mu}$. (a) $\alpha=0$, (b) $\alpha=2$, (c) $\alpha=4$, (d) $\alpha=8$, (e) $\alpha=16$, and (f) $\alpha=32$.}
  \label{fig:vel}
\end{figure}

All the results are normalized using the steady centerline velocity from the reference solution, $u_{\rm ref}(0,0)$. 
An overall good agreement is observed between the simulation and the reference results (Figure \ref{fig:vel}), particularly at smaller Womersley numbers.
Accurate computation of the velocity profile at larger Womersley numbers requires the use finer grids as sharper gradients are developed in those regimes. 

\begin{figure}
  \centering
\begingroup
  \makeatletter
  \providecommand\color[2][]{%
    \GenericError{(gnuplot) \space\space\space\@spaces}{%
      Package color not loaded in conjunction with
      terminal option `colourtext'%
    }{See the gnuplot documentation for explanation.%
    }{Either use 'blacktext' in gnuplot or load the package
      color.sty in LaTeX.}%
    \renewcommand\color[2][]{}%
  }%
  \providecommand\includegraphics[2][]{%
    \GenericError{(gnuplot) \space\space\space\@spaces}{%
      Package graphicx or graphics not loaded%
    }{See the gnuplot documentation for explanation.%
    }{The gnuplot epslatex terminal needs graphicx.sty or graphics.sty.}%
    \renewcommand\includegraphics[2][]{}%
  }%
  \providecommand\rotatebox[2]{#2}%
  \@ifundefined{ifGPcolor}{%
    \newif\ifGPcolor
    \GPcolortrue
  }{}%
  \@ifundefined{ifGPblacktext}{%
    \newif\ifGPblacktext
    \GPblacktexttrue
  }{}%
  \let\gplgaddtomacro\g@addto@macro
  \gdef\gplbacktext{}%
  \gdef\gplfronttext{}%
  \makeatother
  \ifGPblacktext
    \def\colorrgb#1{}%
    \def\colorgray#1{}%
  \else
    \ifGPcolor
      \def\colorrgb#1{\color[rgb]{#1}}%
      \def\colorgray#1{\color[gray]{#1}}%
      \expandafter\def\csname LTw\endcsname{\color{white}}%
      \expandafter\def\csname LTb\endcsname{\color{black}}%
      \expandafter\def\csname LTa\endcsname{\color{black}}%
      \expandafter\def\csname LT0\endcsname{\color[rgb]{1,0,0}}%
      \expandafter\def\csname LT1\endcsname{\color[rgb]{0,1,0}}%
      \expandafter\def\csname LT2\endcsname{\color[rgb]{0,0,1}}%
      \expandafter\def\csname LT3\endcsname{\color[rgb]{1,0,1}}%
      \expandafter\def\csname LT4\endcsname{\color[rgb]{0,1,1}}%
      \expandafter\def\csname LT5\endcsname{\color[rgb]{1,1,0}}%
      \expandafter\def\csname LT6\endcsname{\color[rgb]{0,0,0}}%
      \expandafter\def\csname LT7\endcsname{\color[rgb]{1,0.3,0}}%
      \expandafter\def\csname LT8\endcsname{\color[rgb]{0.5,0.5,0.5}}%
    \else
      \def\colorrgb#1{\color{black}}%
      \def\colorgray#1{\color[gray]{#1}}%
      \expandafter\def\csname LTw\endcsname{\color{white}}%
      \expandafter\def\csname LTb\endcsname{\color{black}}%
      \expandafter\def\csname LTa\endcsname{\color{black}}%
      \expandafter\def\csname LT0\endcsname{\color{black}}%
      \expandafter\def\csname LT1\endcsname{\color{black}}%
      \expandafter\def\csname LT2\endcsname{\color{black}}%
      \expandafter\def\csname LT3\endcsname{\color{black}}%
      \expandafter\def\csname LT4\endcsname{\color{black}}%
      \expandafter\def\csname LT5\endcsname{\color{black}}%
      \expandafter\def\csname LT6\endcsname{\color{black}}%
      \expandafter\def\csname LT7\endcsname{\color{black}}%
      \expandafter\def\csname LT8\endcsname{\color{black}}%
    \fi
  \fi
    \setlength{\unitlength}{0.0500bp}%
    \ifx\gptboxheight\undefined%
      \newlength{\gptboxheight}%
      \newlength{\gptboxwidth}%
      \newsavebox{\gptboxtext}%
    \fi%
    \setlength{\fboxrule}{0.5pt}%
    \setlength{\fboxsep}{1pt}%
\begin{picture}(5040.00,3600.00)%
    \gplgaddtomacro\gplbacktext{%
      \csname LTb\endcsname
      \put(1210,550){\makebox(0,0)[r]{\strut{}$0.0001$}}%
      \put(1210,1257){\makebox(0,0)[r]{\strut{}$0.001$}}%
      \put(1210,1965){\makebox(0,0)[r]{\strut{}$0.01$}}%
      \put(1210,2672){\makebox(0,0)[r]{\strut{}$0.1$}}%
      \put(1210,3379){\makebox(0,0)[r]{\strut{}$1$}}%
      \put(1709,330){\makebox(0,0){\strut{}$2$}}%
      \put(2442,330){\makebox(0,0){\strut{}$4$}}%
      \put(3176,330){\makebox(0,0){\strut{}$8$}}%
      \put(3909,330){\makebox(0,0){\strut{}$16$}}%
      \put(4643,330){\makebox(0,0){\strut{}$32$}}%
      \put(3653,1116){\makebox(0,0)[l]{\strut{}$q_r^*$}}%
      \put(3653,2672){\makebox(0,0)[l]{\strut{}$-q_i^*$}}%
    }%
    \gplgaddtomacro\gplfronttext{%
      \csname LTb\endcsname
      \put(473,1964){\rotatebox{-270}{\makebox(0,0){\strut{}$q_r^*;\; -q_i^*$}}}%
      \put(2992,154){\makebox(0,0){\strut{}$\alpha$}}%
      \put(2002,1383){\makebox(0,0)[r]{\strut{}Ref.}}%
      \csname LTb\endcsname
      \put(2002,1163){\makebox(0,0)[r]{\strut{}M1}}%
      \csname LTb\endcsname
      \put(2002,943){\makebox(0,0)[r]{\strut{}M2}}%
      \csname LTb\endcsname
      \put(2002,723){\makebox(0,0)[r]{\strut{}M3}}%
    }%
    \gplbacktext
    \put(0,0){\includegraphics[width={252.00bp},height={180.00bp}]{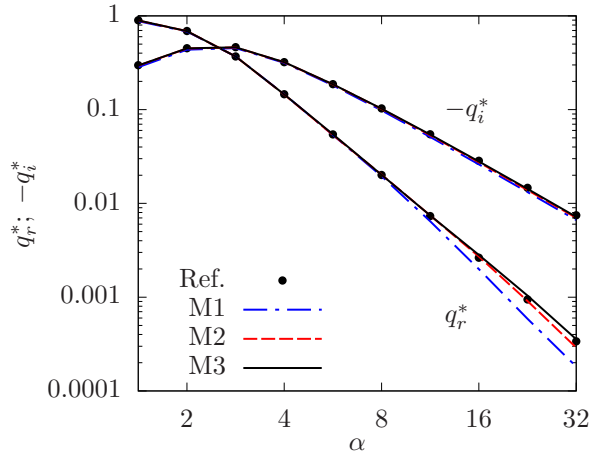}}%
    \gplfronttext
  \end{picture}%
\endgroup
   \caption{The normalized real ($q_r^*$) and imaginary ($q_i^*$) flow rate versus Womersley number $\alpha=R\sqrt{\rho\omega/\mu}$ for three meshes, M1 (dashed-dot), M2 (dashed), and M3 (solid) and their comparison against the reference results (dots) from Eq. \eqref{wom_q}.}
  \label{fig:flow}    
\end{figure}

The real and imaginary components of the predicted flow rate are compared against the analytical solution of Womersley that is 
\begin{equation}
   q_{\rm ref}(\alpha) = \left\{ \begin{array}{lr} 
    \displaystyle \frac{h\pi R^4}{8\mu L}, &  \alpha = 0, \\ 
    & \\
       \displaystyle \frac{-\hat{j}\pi h R^4}{L\mu \alpha^2} \left[1+\frac{2\hat j^{\frac1 2} J_1(\hat j^{\frac3 2}\alpha)}{\alpha J_0(\hat j^{\frac3 2}\alpha)}\right],  & \alpha \neq 0,
    \end{array}\right.
    \label{wom_q}
\end{equation}
in Figure \ref{fig:flow}.
The results, which are normalized based on the reference steady flow rate ($q_{ref}(0)$), show a good agreement with the reference results. 
As correctly computed, the portion of the flow rate that is in phase with pressure drops faster at larger $\alpha$ than the out of phase component (i.e., $q_r \to 0$ faster than $q_i \to 0$ as $\alpha \to \infty$). 
Thus, the fact that pressure and flow are in phase and out of phase by $\pi/2$ at small and large $\alpha$, respectively, are correctly captured by the proposed formulation. 

\begin{figure}
  \centering
\begingroup
  \makeatletter
  \providecommand\color[2][]{%
    \GenericError{(gnuplot) \space\space\space\@spaces}{%
      Package color not loaded in conjunction with
      terminal option `colourtext'%
    }{See the gnuplot documentation for explanation.%
    }{Either use 'blacktext' in gnuplot or load the package
      color.sty in LaTeX.}%
    \renewcommand\color[2][]{}%
  }%
  \providecommand\includegraphics[2][]{%
    \GenericError{(gnuplot) \space\space\space\@spaces}{%
      Package graphicx or graphics not loaded%
    }{See the gnuplot documentation for explanation.%
    }{The gnuplot epslatex terminal needs graphicx.sty or graphics.sty.}%
    \renewcommand\includegraphics[2][]{}%
  }%
  \providecommand\rotatebox[2]{#2}%
  \@ifundefined{ifGPcolor}{%
    \newif\ifGPcolor
    \GPcolortrue
  }{}%
  \@ifundefined{ifGPblacktext}{%
    \newif\ifGPblacktext
    \GPblacktexttrue
  }{}%
  \let\gplgaddtomacro\g@addto@macro
  \gdef\gplbacktext{}%
  \gdef\gplfronttext{}%
  \makeatother
  \ifGPblacktext
    \def\colorrgb#1{}%
    \def\colorgray#1{}%
  \else
    \ifGPcolor
      \def\colorrgb#1{\color[rgb]{#1}}%
      \def\colorgray#1{\color[gray]{#1}}%
      \expandafter\def\csname LTw\endcsname{\color{white}}%
      \expandafter\def\csname LTb\endcsname{\color{black}}%
      \expandafter\def\csname LTa\endcsname{\color{black}}%
      \expandafter\def\csname LT0\endcsname{\color[rgb]{1,0,0}}%
      \expandafter\def\csname LT1\endcsname{\color[rgb]{0,1,0}}%
      \expandafter\def\csname LT2\endcsname{\color[rgb]{0,0,1}}%
      \expandafter\def\csname LT3\endcsname{\color[rgb]{1,0,1}}%
      \expandafter\def\csname LT4\endcsname{\color[rgb]{0,1,1}}%
      \expandafter\def\csname LT5\endcsname{\color[rgb]{1,1,0}}%
      \expandafter\def\csname LT6\endcsname{\color[rgb]{0,0,0}}%
      \expandafter\def\csname LT7\endcsname{\color[rgb]{1,0.3,0}}%
      \expandafter\def\csname LT8\endcsname{\color[rgb]{0.5,0.5,0.5}}%
    \else
      \def\colorrgb#1{\color{black}}%
      \def\colorgray#1{\color[gray]{#1}}%
      \expandafter\def\csname LTw\endcsname{\color{white}}%
      \expandafter\def\csname LTb\endcsname{\color{black}}%
      \expandafter\def\csname LTa\endcsname{\color{black}}%
      \expandafter\def\csname LT0\endcsname{\color{black}}%
      \expandafter\def\csname LT1\endcsname{\color{black}}%
      \expandafter\def\csname LT2\endcsname{\color{black}}%
      \expandafter\def\csname LT3\endcsname{\color{black}}%
      \expandafter\def\csname LT4\endcsname{\color{black}}%
      \expandafter\def\csname LT5\endcsname{\color{black}}%
      \expandafter\def\csname LT6\endcsname{\color{black}}%
      \expandafter\def\csname LT7\endcsname{\color{black}}%
      \expandafter\def\csname LT8\endcsname{\color{black}}%
    \fi
  \fi
    \setlength{\unitlength}{0.0500bp}%
    \ifx\gptboxheight\undefined%
      \newlength{\gptboxheight}%
      \newlength{\gptboxwidth}%
      \newsavebox{\gptboxtext}%
    \fi%
    \setlength{\fboxrule}{0.5pt}%
    \setlength{\fboxsep}{1pt}%
\begin{picture}(9360.00,3600.00)%
    \gplgaddtomacro\gplbacktext{%
      \csname LTb\endcsname
      \put(991,1344){\makebox(0,0)[r]{\strut{}$0.01$}}%
      \put(991,2752){\makebox(0,0)[r]{\strut{}$0.1$}}%
      \put(1518,140){\makebox(0,0){\strut{}$2$}}%
      \put(2308,140){\makebox(0,0){\strut{}$4$}}%
      \put(3099,140){\makebox(0,0){\strut{}$8$}}%
      \put(3889,140){\makebox(0,0){\strut{}$16$}}%
      \put(4679,140){\makebox(0,0){\strut{}$32$}}%
      \put(4323,684){\makebox(0,0)[l]{\strut{}(a)}}%
    }%
    \gplgaddtomacro\gplfronttext{%
      \csname LTb\endcsname
      \put(386,1979){\rotatebox{-270}{\makebox(0,0){\strut{}$\|u-u_{\mathrm ref}\|/\|u_{\mathrm ref}\|$}}}%
      \put(2901,-36){\makebox(0,0){\strut{}$\alpha$}}%
      \csname LTb\endcsname
      \put(1519,3426){\makebox(0,0)[r]{\strut{}M1}}%
      \csname LTb\endcsname
      \put(1519,3206){\makebox(0,0)[r]{\strut{}M2}}%
      \csname LTb\endcsname
      \put(1519,2986){\makebox(0,0)[r]{\strut{}M3}}%
    }%
    \gplgaddtomacro\gplbacktext{%
      \csname LTb\endcsname
      \put(5671,360){\makebox(0,0)[r]{\strut{}$0$}}%
      \put(5671,1170){\makebox(0,0)[r]{\strut{}$500$}}%
      \put(5671,1980){\makebox(0,0)[r]{\strut{}$1000$}}%
      \put(5671,2789){\makebox(0,0)[r]{\strut{}$1500$}}%
      \put(5671,3599){\makebox(0,0)[r]{\strut{}$2000$}}%
      \put(6198,140){\makebox(0,0){\strut{}$2$}}%
      \put(6988,140){\makebox(0,0){\strut{}$4$}}%
      \put(7779,140){\makebox(0,0){\strut{}$8$}}%
      \put(8569,140){\makebox(0,0){\strut{}$16$}}%
      \put(9359,140){\makebox(0,0){\strut{}$32$}}%
      \put(9003,684){\makebox(0,0)[l]{\strut{}(b)}}%
    }%
    \gplgaddtomacro\gplfronttext{%
      \csname LTb\endcsname
      \put(5066,1979){\rotatebox{-270}{\makebox(0,0){\strut{}$N_{\mathrm itr}$}}}%
      \put(7581,-36){\makebox(0,0){\strut{}$\alpha$}}%
    }%
    \gplbacktext
    \put(0,0){\includegraphics[width={468.00bp},height={180.00bp}]{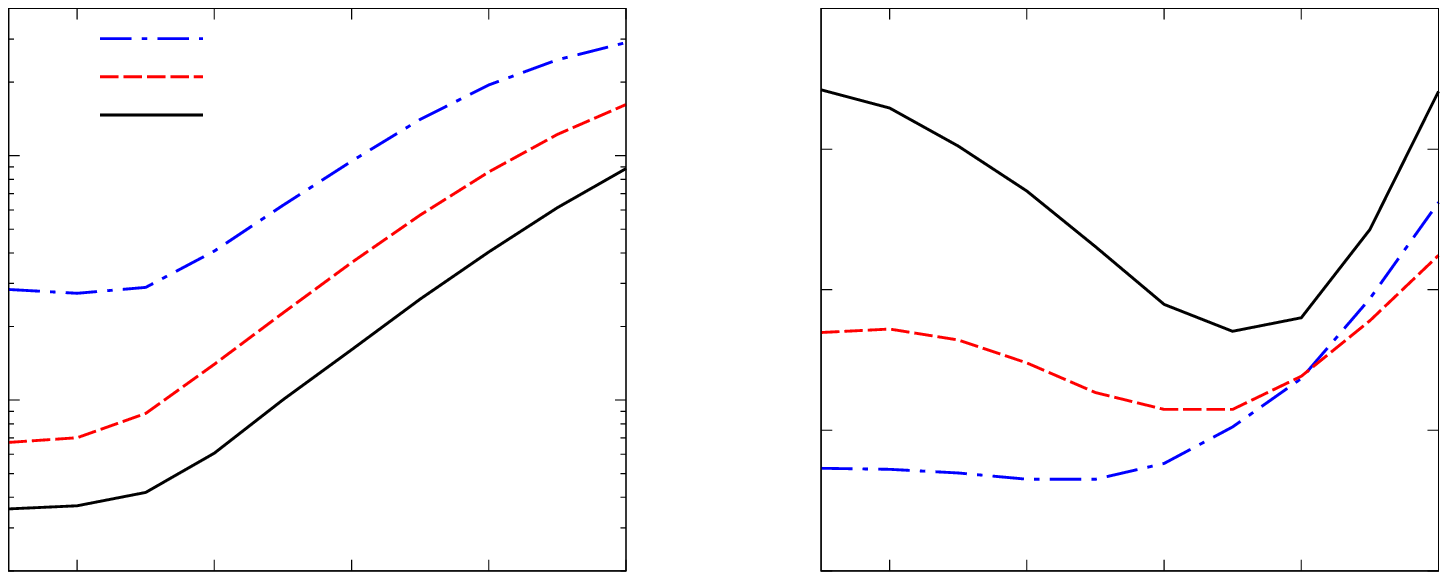}}%
    \gplfronttext
  \end{picture}%
\endgroup
   \caption{The error (a) and the number of conjugate gradient iterations (b) used for simulating the oscillating pipe flow versus Womersley number $\alpha=R\sqrt{\rho\omega/\mu}$ for three meshes, M1 (dashed-dot), M2 (dashed), and M3 (solid).}
  \label{fig:alpha}    
\end{figure}

The change in error as a function of Womersley number (or oscillation frequency $\omega$) and mesh resolution is captured more explicitly in Figure \ref{fig:alpha}(a). 
Similar to what was observed earlier, increasing the Womersley number or the element size increases the error in the solution. 
The increase in error with $\alpha$ only occurs beyond a certain threshold since the curvature in the fluid velocity field remains relatively unchanged at small $\alpha$.

To better understand the computational cost of the proposed method, the number of iterations of the linear solver $N_{\rm itr}$ is plotted for these computations in Figure \ref{fig:alpha}(b). 
These results show that $N_{\rm itr} = \mathcal O(10^3)$ as $\alpha$ is increased from zero to 32, indicating the relative independence of the computational cost as the flow becomes more oscillatory.
This independence is in contrast to the previous Bubnov-Galerkin formulation in \cite{meng2020time} that showed an increase in cost at high Womersley numbers due to the ill-conditioning of the tangent matrix owing to its zero diagonal block. 
The higher number of iterations required for the finer grid in Figure \ref{fig:alpha}-(b) is a consequence of the broader range of eigenvalues in the tangent matrix. 
Adopting a multigrid preconditioner in future studies should break the dependence of $N_{\rm itr}$ on $L/H$, thus leading to an overall cost that linearly scales with the number of nodes in the mesh. 

\begin{figure}
  \centering
\begingroup
  \makeatletter
  \providecommand\color[2][]{%
    \GenericError{(gnuplot) \space\space\space\@spaces}{%
      Package color not loaded in conjunction with
      terminal option `colourtext'%
    }{See the gnuplot documentation for explanation.%
    }{Either use 'blacktext' in gnuplot or load the package
      color.sty in LaTeX.}%
    \renewcommand\color[2][]{}%
  }%
  \providecommand\includegraphics[2][]{%
    \GenericError{(gnuplot) \space\space\space\@spaces}{%
      Package graphicx or graphics not loaded%
    }{See the gnuplot documentation for explanation.%
    }{The gnuplot epslatex terminal needs graphicx.sty or graphics.sty.}%
    \renewcommand\includegraphics[2][]{}%
  }%
  \providecommand\rotatebox[2]{#2}%
  \@ifundefined{ifGPcolor}{%
    \newif\ifGPcolor
    \GPcolortrue
  }{}%
  \@ifundefined{ifGPblacktext}{%
    \newif\ifGPblacktext
    \GPblacktexttrue
  }{}%
  \let\gplgaddtomacro\g@addto@macro
  \gdef\gplbacktext{}%
  \gdef\gplfronttext{}%
  \makeatother
  \ifGPblacktext
    \def\colorrgb#1{}%
    \def\colorgray#1{}%
  \else
    \ifGPcolor
      \def\colorrgb#1{\color[rgb]{#1}}%
      \def\colorgray#1{\color[gray]{#1}}%
      \expandafter\def\csname LTw\endcsname{\color{white}}%
      \expandafter\def\csname LTb\endcsname{\color{black}}%
      \expandafter\def\csname LTa\endcsname{\color{black}}%
      \expandafter\def\csname LT0\endcsname{\color[rgb]{1,0,0}}%
      \expandafter\def\csname LT1\endcsname{\color[rgb]{0,1,0}}%
      \expandafter\def\csname LT2\endcsname{\color[rgb]{0,0,1}}%
      \expandafter\def\csname LT3\endcsname{\color[rgb]{1,0,1}}%
      \expandafter\def\csname LT4\endcsname{\color[rgb]{0,1,1}}%
      \expandafter\def\csname LT5\endcsname{\color[rgb]{1,1,0}}%
      \expandafter\def\csname LT6\endcsname{\color[rgb]{0,0,0}}%
      \expandafter\def\csname LT7\endcsname{\color[rgb]{1,0.3,0}}%
      \expandafter\def\csname LT8\endcsname{\color[rgb]{0.5,0.5,0.5}}%
    \else
      \def\colorrgb#1{\color{black}}%
      \def\colorgray#1{\color[gray]{#1}}%
      \expandafter\def\csname LTw\endcsname{\color{white}}%
      \expandafter\def\csname LTb\endcsname{\color{black}}%
      \expandafter\def\csname LTa\endcsname{\color{black}}%
      \expandafter\def\csname LT0\endcsname{\color{black}}%
      \expandafter\def\csname LT1\endcsname{\color{black}}%
      \expandafter\def\csname LT2\endcsname{\color{black}}%
      \expandafter\def\csname LT3\endcsname{\color{black}}%
      \expandafter\def\csname LT4\endcsname{\color{black}}%
      \expandafter\def\csname LT5\endcsname{\color{black}}%
      \expandafter\def\csname LT6\endcsname{\color{black}}%
      \expandafter\def\csname LT7\endcsname{\color{black}}%
      \expandafter\def\csname LT8\endcsname{\color{black}}%
    \fi
  \fi
    \setlength{\unitlength}{0.0500bp}%
    \ifx\gptboxheight\undefined%
      \newlength{\gptboxheight}%
      \newlength{\gptboxwidth}%
      \newsavebox{\gptboxtext}%
    \fi%
    \setlength{\fboxrule}{0.5pt}%
    \setlength{\fboxsep}{1pt}%
\begin{picture}(5040.00,3600.00)%
    \gplgaddtomacro\gplbacktext{%
      \csname LTb\endcsname
      \put(1078,660){\makebox(0,0)[r]{\strut{}1e-06}}%
      \put(1078,1113){\makebox(0,0)[r]{\strut{}1e-04}}%
      \put(1078,1566){\makebox(0,0)[r]{\strut{}1e-02}}%
      \put(1078,2020){\makebox(0,0)[r]{\strut{}1e+00}}%
      \put(1078,2473){\makebox(0,0)[r]{\strut{}1e+02}}%
      \put(1078,2926){\makebox(0,0)[r]{\strut{}1e+04}}%
      \put(1078,3379){\makebox(0,0)[r]{\strut{}1e+06}}%
      \put(1210,440){\makebox(0,0){\strut{}$0$}}%
      \put(1782,440){\makebox(0,0){\strut{}$500$}}%
      \put(2354,440){\makebox(0,0){\strut{}$1000$}}%
      \put(2927,440){\makebox(0,0){\strut{}$1500$}}%
      \put(3499,440){\makebox(0,0){\strut{}$2000$}}%
      \put(4071,440){\makebox(0,0){\strut{}$2500$}}%
      \put(4643,440){\makebox(0,0){\strut{}$3000$}}%
    }%
    \gplgaddtomacro\gplfronttext{%
      \csname LTb\endcsname
      \put(341,2019){\rotatebox{-270}{\makebox(0,0){\strut{}$\|R\|^2$}}}%
      \put(2926,154){\makebox(0,0){\strut{}$i$}}%
      \put(3656,3206){\makebox(0,0)[r]{\strut{}$\alpha=0$}}%
      \csname LTb\endcsname
      \put(3656,2986){\makebox(0,0)[r]{\strut{}$\alpha=5$}}%
      \csname LTb\endcsname
      \put(3656,2766){\makebox(0,0)[r]{\strut{}$\alpha=10$}}%
    }%
    \gplbacktext
    \put(0,0){\includegraphics[width={252.00bp},height={180.00bp}]{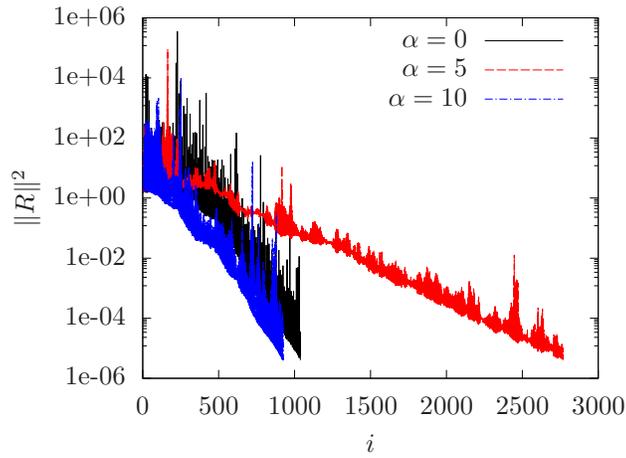}}%
    \gplfronttext
  \end{picture}%
\endgroup
   \caption{The error (a) and the number of conjugate gradient iterations (b) used for simulating the oscillating pipe flow versus the stabilization parameter $c$ in Eq. \eqref{tau}. The results are obtained at $\alpha=1$ for three meshes, M1 (dashed-dot), M2 (dashed), and M3 (solid).
   Similar behavior is observed at other values of $\alpha$.}
  \label{fig:ci}
\end{figure}

As stated earlier, the proposed stabilized method has one adjustable parameter $c$ that appeared in Eq. \eqref{tau}. 
One of the attractive properties of the proposed scheme is that its accuracy hardly depends on this parameter. 
As shown in Figure \ref{fig:ci}, changing $c$ by orders of magnitude has little effect on the error.
However, the number of linear solver iterations $N_{\rm itr}$ and thus the overall cost do depend on $c$. 
These results, which are obtained at $\alpha=1$, are relatively independent of $\alpha$.
Therefore, it is expected that the proposed $c=2^{-5}$ to produce accurate results across the board even if it produces a suboptimal convergence rate when an element type other than those tested here (i.e., bilateral and triangular) are used. 

\begin{figure}
  \centering
\begingroup
  \makeatletter
  \providecommand\color[2][]{%
    \GenericError{(gnuplot) \space\space\space\@spaces}{%
      Package color not loaded in conjunction with
      terminal option `colourtext'%
    }{See the gnuplot documentation for explanation.%
    }{Either use 'blacktext' in gnuplot or load the package
      color.sty in LaTeX.}%
    \renewcommand\color[2][]{}%
  }%
  \providecommand\includegraphics[2][]{%
    \GenericError{(gnuplot) \space\space\space\@spaces}{%
      Package graphicx or graphics not loaded%
    }{See the gnuplot documentation for explanation.%
    }{The gnuplot epslatex terminal needs graphicx.sty or graphics.sty.}%
    \renewcommand\includegraphics[2][]{}%
  }%
  \providecommand\rotatebox[2]{#2}%
  \@ifundefined{ifGPcolor}{%
    \newif\ifGPcolor
    \GPcolortrue
  }{}%
  \@ifundefined{ifGPblacktext}{%
    \newif\ifGPblacktext
    \GPblacktexttrue
  }{}%
  \let\gplgaddtomacro\g@addto@macro
  \gdef\gplbacktext{}%
  \gdef\gplfronttext{}%
  \makeatother
  \ifGPblacktext
    \def\colorrgb#1{}%
    \def\colorgray#1{}%
  \else
    \ifGPcolor
      \def\colorrgb#1{\color[rgb]{#1}}%
      \def\colorgray#1{\color[gray]{#1}}%
      \expandafter\def\csname LTw\endcsname{\color{white}}%
      \expandafter\def\csname LTb\endcsname{\color{black}}%
      \expandafter\def\csname LTa\endcsname{\color{black}}%
      \expandafter\def\csname LT0\endcsname{\color[rgb]{1,0,0}}%
      \expandafter\def\csname LT1\endcsname{\color[rgb]{0,1,0}}%
      \expandafter\def\csname LT2\endcsname{\color[rgb]{0,0,1}}%
      \expandafter\def\csname LT3\endcsname{\color[rgb]{1,0,1}}%
      \expandafter\def\csname LT4\endcsname{\color[rgb]{0,1,1}}%
      \expandafter\def\csname LT5\endcsname{\color[rgb]{1,1,0}}%
      \expandafter\def\csname LT6\endcsname{\color[rgb]{0,0,0}}%
      \expandafter\def\csname LT7\endcsname{\color[rgb]{1,0.3,0}}%
      \expandafter\def\csname LT8\endcsname{\color[rgb]{0.5,0.5,0.5}}%
    \else
      \def\colorrgb#1{\color{black}}%
      \def\colorgray#1{\color[gray]{#1}}%
      \expandafter\def\csname LTw\endcsname{\color{white}}%
      \expandafter\def\csname LTb\endcsname{\color{black}}%
      \expandafter\def\csname LTa\endcsname{\color{black}}%
      \expandafter\def\csname LT0\endcsname{\color{black}}%
      \expandafter\def\csname LT1\endcsname{\color{black}}%
      \expandafter\def\csname LT2\endcsname{\color{black}}%
      \expandafter\def\csname LT3\endcsname{\color{black}}%
      \expandafter\def\csname LT4\endcsname{\color{black}}%
      \expandafter\def\csname LT5\endcsname{\color{black}}%
      \expandafter\def\csname LT6\endcsname{\color{black}}%
      \expandafter\def\csname LT7\endcsname{\color{black}}%
      \expandafter\def\csname LT8\endcsname{\color{black}}%
    \fi
  \fi
    \setlength{\unitlength}{0.0500bp}%
    \ifx\gptboxheight\undefined%
      \newlength{\gptboxheight}%
      \newlength{\gptboxwidth}%
      \newsavebox{\gptboxtext}%
    \fi%
    \setlength{\fboxrule}{0.5pt}%
    \setlength{\fboxsep}{1pt}%
\begin{picture}(5040.00,3600.00)%
    \gplgaddtomacro\gplbacktext{%
      \csname LTb\endcsname
      \put(814,834){\makebox(0,0)[r]{\strut{}$0.1$}}%
      \put(814,1777){\makebox(0,0)[r]{\strut{}$1$}}%
      \put(814,2720){\makebox(0,0)[r]{\strut{}$10$}}%
      \put(1357,330){\makebox(0,0){\strut{}$2$}}%
      \put(2178,330){\makebox(0,0){\strut{}$4$}}%
      \put(3000,330){\makebox(0,0){\strut{}$8$}}%
      \put(3821,330){\makebox(0,0){\strut{}$16$}}%
      \put(4643,330){\makebox(0,0){\strut{}$32$}}%
    }%
    \gplgaddtomacro\gplfronttext{%
      \csname LTb\endcsname
      \put(341,1964){\rotatebox{-270}{\makebox(0,0){\strut{}Wall clock time (sec)}}}%
      \put(2794,154){\makebox(0,0){\strut{}$\alpha$}}%
      \csname LTb\endcsname
      \put(3656,1163){\makebox(0,0)[r]{\strut{}M1}}%
      \csname LTb\endcsname
      \put(3656,943){\makebox(0,0)[r]{\strut{}M2}}%
      \csname LTb\endcsname
      \put(3656,723){\makebox(0,0)[r]{\strut{}M3}}%
    }%
    \gplbacktext
    \put(0,0){\includegraphics[width={252.00bp},height={180.00bp}]{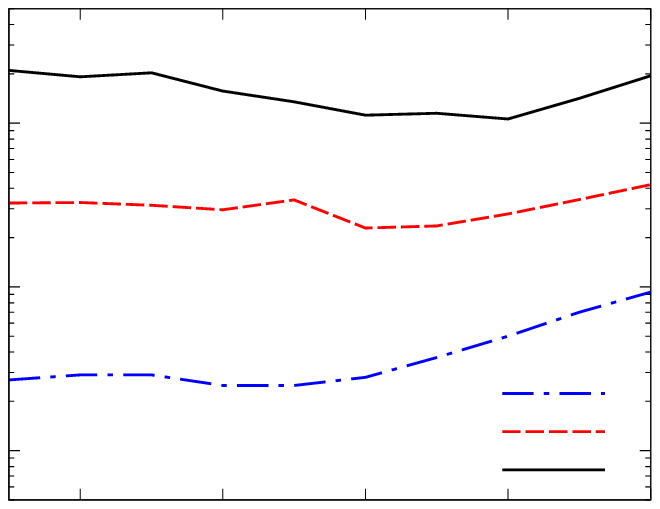}}%
    \gplfronttext
  \end{picture}%
\endgroup
   \caption{The simulation turnover time in second as a function of Womersley number $\alpha=R\sqrt{\rho\omega/\mu}$ for three meshes, M1 (dashed-dot), M2 (dashed), and M3 (solid).}
  \label{fig:cost}    
\end{figure}

The computational performance of the proposed scheme, both in terms of CPU-hours and solution turnover time, is far shorter than a standard CFD solver and the previously proposed Bubnov-Galerkin formulation. 
Taking mesh M2 at $\alpha^2 = 8\pi$ as an example, the present formulation takes 3 seconds to produce a solution using 32 processors, which is roughly 0.027 CPU-hours (Figure \ref{fig:cost}). 
Similar simulation performed using the Bubnov-Galerkin formulation and a standard CFD solver will require $1.25$ and over $300$ CPU-hours, respectively \cite{meng2020time}. 
The four orders of magnitude reduction in cost in comparison to the standard CFD solver is roughly the same as the ratio between the number of time steps for the CFD solver and the number of solution modes for the present formulation, which is 10,000 to 1. 
Even if the number of modes required for an accurate representation of the boundary conditions is $\mathcal O(10)$, the present formulation will reduce the cost by three orders of magnitude. 
The solution turnover time of this case, on the other hand, will always be four orders of magnitude shorter than that of the conventional CFD solver given that all the modes can be computed in parallel for the linear Stokes equations. 

\subsection{Patient-specific Glenn geometry} 
\begin{figure}
  \centering
  \includegraphics[width=1\textwidth]{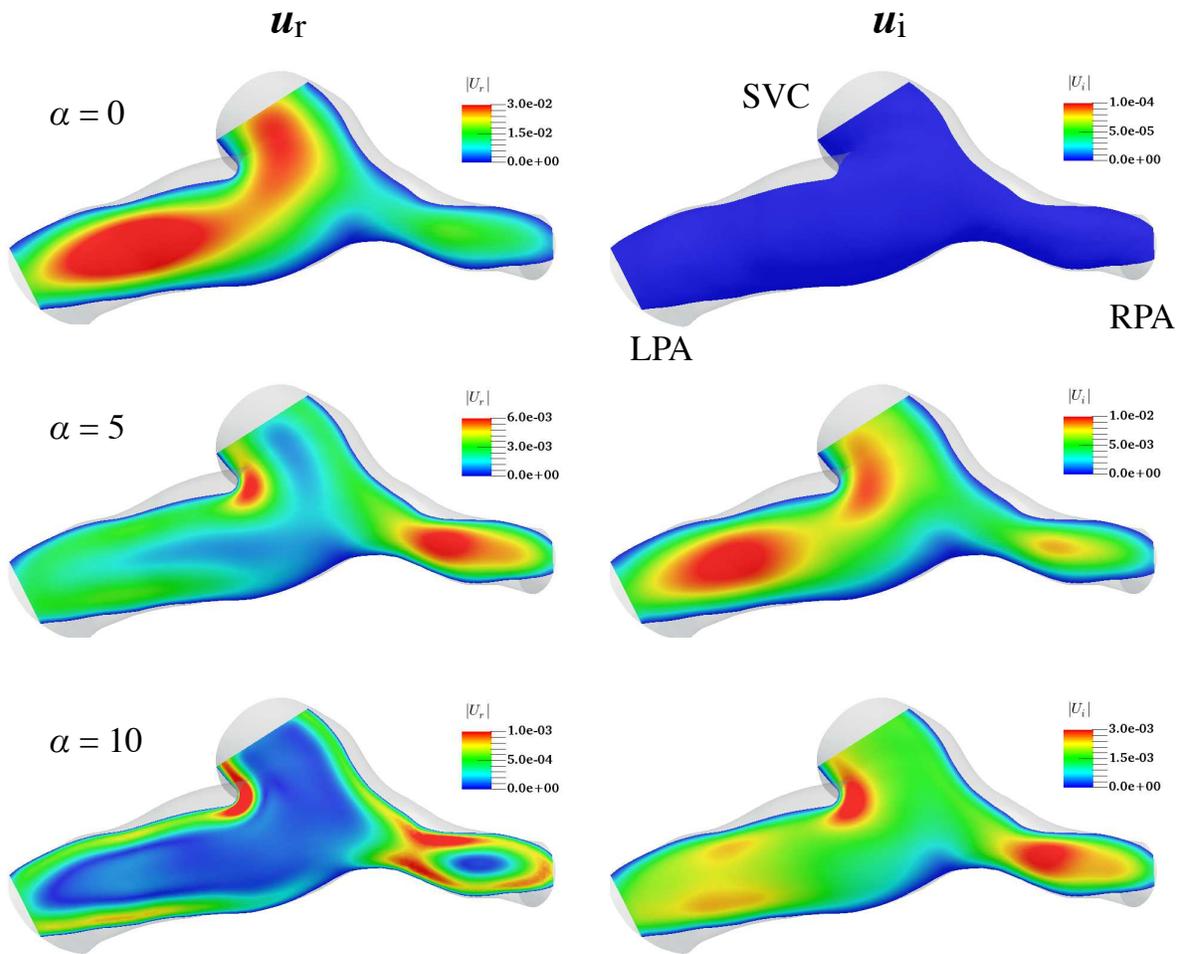}
   \caption{The real (left column) and imaginary (right column) components of velocity shown on a 2D slice of the Glenn geometry. The first row corresponds to the steady flow ($\alpha=0$), the second row corresponds to moderately oscillatory flow ($\alpha=5$), and the last row corresponds to a highly oscillatory flow ($\alpha=10$). }
  \label{fig:stubby}    
\end{figure}

To demonstrate the applicability of the proposed method to more complex problems, a patient-specific geometry obtained from a patient undergoing Glenn operation is considered \cite{arbia2014numerical}.
There are three inflow-outflow boundaries, namely superior vena cava (SVC) and the left and right pulmonary arteries (LPA and RPA).
Zero Neumann boundary condition is imposed at the LPA and RPA and a non-slip boundary condition is imposed at the walls (Figure \ref{fig:stubby}).
To derive the flow, a unit oscillatory traction is imposed at the SVC with a frequency that is adjusted to obtain different Womersley numbers. 
The Womersley number, in this case, is defined as $\alpha = D_{\rm h}\sqrt{\rho\omega/\mu}$, where $D_{\rm h}$ is the hydraulic diameter of the SVC. 
The geometry is discretized using 988,747 linear tetrahedral elements, resulting in 163,791 nodes. 

The significant dependency of the solution to the Womersley number is shown in Figure \ref{fig:stubby} for the simulations performed at $\alpha=0,$ 5, and 10.
Despite geometrical complexity, the overall variation in the real and imaginary component of velocity follows that of the pipe flow case. 
For the steady case with $\alpha=0$, $\bl u_i =\bl 0$ and $\bl u_r$ forms a parabolic profile at the annular cross sections. 
As $\alpha$ increases, $\bl u_r$ develops peaks in the near wall region while $\bl u_i$ grows in relative magnitude to create a bulk flow at the annular cross sections. 

\begin{figure}
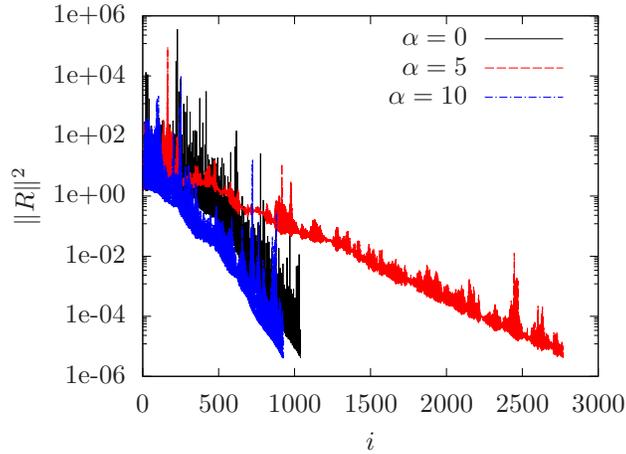

  \centering
\begingroup
  \makeatletter
  \providecommand\color[2][]{%
    \GenericError{(gnuplot) \space\space\space\@spaces}{%
      Package color not loaded in conjunction with
      terminal option `colourtext'%
    }{See the gnuplot documentation for explanation.%
    }{Either use 'blacktext' in gnuplot or load the package
      color.sty in LaTeX.}%
    \renewcommand\color[2][]{}%
  }%
  \providecommand\includegraphics[2][]{%
    \GenericError{(gnuplot) \space\space\space\@spaces}{%
      Package graphicx or graphics not loaded%
    }{See the gnuplot documentation for explanation.%
    }{The gnuplot epslatex terminal needs graphicx.sty or graphics.sty.}%
    \renewcommand\includegraphics[2][]{}%
  }%
  \providecommand\rotatebox[2]{#2}%
  \@ifundefined{ifGPcolor}{%
    \newif\ifGPcolor
    \GPcolortrue
  }{}%
  \@ifundefined{ifGPblacktext}{%
    \newif\ifGPblacktext
    \GPblacktexttrue
  }{}%
  \let\gplgaddtomacro\g@addto@macro
  \gdef\gplbacktext{}%
  \gdef\gplfronttext{}%
  \makeatother
  \ifGPblacktext
    \def\colorrgb#1{}%
    \def\colorgray#1{}%
  \else
    \ifGPcolor
      \def\colorrgb#1{\color[rgb]{#1}}%
      \def\colorgray#1{\color[gray]{#1}}%
      \expandafter\def\csname LTw\endcsname{\color{white}}%
      \expandafter\def\csname LTb\endcsname{\color{black}}%
      \expandafter\def\csname LTa\endcsname{\color{black}}%
      \expandafter\def\csname LT0\endcsname{\color[rgb]{1,0,0}}%
      \expandafter\def\csname LT1\endcsname{\color[rgb]{0,1,0}}%
      \expandafter\def\csname LT2\endcsname{\color[rgb]{0,0,1}}%
      \expandafter\def\csname LT3\endcsname{\color[rgb]{1,0,1}}%
      \expandafter\def\csname LT4\endcsname{\color[rgb]{0,1,1}}%
      \expandafter\def\csname LT5\endcsname{\color[rgb]{1,1,0}}%
      \expandafter\def\csname LT6\endcsname{\color[rgb]{0,0,0}}%
      \expandafter\def\csname LT7\endcsname{\color[rgb]{1,0.3,0}}%
      \expandafter\def\csname LT8\endcsname{\color[rgb]{0.5,0.5,0.5}}%
    \else
      \def\colorrgb#1{\color{black}}%
      \def\colorgray#1{\color[gray]{#1}}%
      \expandafter\def\csname LTw\endcsname{\color{white}}%
      \expandafter\def\csname LTb\endcsname{\color{black}}%
      \expandafter\def\csname LTa\endcsname{\color{black}}%
      \expandafter\def\csname LT0\endcsname{\color{black}}%
      \expandafter\def\csname LT1\endcsname{\color{black}}%
      \expandafter\def\csname LT2\endcsname{\color{black}}%
      \expandafter\def\csname LT3\endcsname{\color{black}}%
      \expandafter\def\csname LT4\endcsname{\color{black}}%
      \expandafter\def\csname LT5\endcsname{\color{black}}%
      \expandafter\def\csname LT6\endcsname{\color{black}}%
      \expandafter\def\csname LT7\endcsname{\color{black}}%
      \expandafter\def\csname LT8\endcsname{\color{black}}%
    \fi
  \fi
    \setlength{\unitlength}{0.0500bp}%
    \ifx\gptboxheight\undefined%
      \newlength{\gptboxheight}%
      \newlength{\gptboxwidth}%
      \newsavebox{\gptboxtext}%
    \fi%
    \setlength{\fboxrule}{0.5pt}%
    \setlength{\fboxsep}{1pt}%
\begin{picture}(5040.00,3600.00)%
    \gplgaddtomacro\gplbacktext{%
      \csname LTb\endcsname
      \put(1078,660){\makebox(0,0)[r]{\strut{}1e-06}}%
      \put(1078,1113){\makebox(0,0)[r]{\strut{}1e-04}}%
      \put(1078,1566){\makebox(0,0)[r]{\strut{}1e-02}}%
      \put(1078,2020){\makebox(0,0)[r]{\strut{}1e+00}}%
      \put(1078,2473){\makebox(0,0)[r]{\strut{}1e+02}}%
      \put(1078,2926){\makebox(0,0)[r]{\strut{}1e+04}}%
      \put(1078,3379){\makebox(0,0)[r]{\strut{}1e+06}}%
      \put(1210,440){\makebox(0,0){\strut{}$0$}}%
      \put(1782,440){\makebox(0,0){\strut{}$500$}}%
      \put(2354,440){\makebox(0,0){\strut{}$1000$}}%
      \put(2927,440){\makebox(0,0){\strut{}$1500$}}%
      \put(3499,440){\makebox(0,0){\strut{}$2000$}}%
      \put(4071,440){\makebox(0,0){\strut{}$2500$}}%
      \put(4643,440){\makebox(0,0){\strut{}$3000$}}%
    }%
    \gplgaddtomacro\gplfronttext{%
      \csname LTb\endcsname
      \put(341,2019){\rotatebox{-270}{\makebox(0,0){\strut{}$\|R\|^2$}}}%
      \put(2926,154){\makebox(0,0){\strut{}$i$}}%
      \put(3656,3206){\makebox(0,0)[r]{\strut{}$\alpha=0$}}%
      \csname LTb\endcsname
      \put(3656,2986){\makebox(0,0)[r]{\strut{}$\alpha=5$}}%
      \csname LTb\endcsname
      \put(3656,2766){\makebox(0,0)[r]{\strut{}$\alpha=10$}}%
    }%
    \gplbacktext
    \put(0,0){\includegraphics[width={252.00bp},height={180.00bp}]{fcg}}%
    \gplfronttext
  \end{picture}%
\endgroup
  \caption{The norm of the residual in the conjugate gradient solver as a function of the iteration number. The linear system is obtained from the Glenn geometry at three different Womersley numbers  (solid; $\alpha=0$, dashed; $\alpha=5$, and dashed-dot; $\alpha=10$).}
  \label{fig:cg}
\end{figure}

The convergence behavior of the conjugate gradient solver for the three simulations at $\alpha=0$, 5, and 10 are shown in Figure \ref{fig:cg}. 
Even though there are intermittent increases in the error, all cases show an exponential convergence rate with $\alpha=5$ case requiring the largest number of iterations. 

\begin{figure}
  \centering
\begingroup
  \makeatletter
  \providecommand\color[2][]{%
    \GenericError{(gnuplot) \space\space\space\@spaces}{%
      Package color not loaded in conjunction with
      terminal option `colourtext'%
    }{See the gnuplot documentation for explanation.%
    }{Either use 'blacktext' in gnuplot or load the package
      color.sty in LaTeX.}%
    \renewcommand\color[2][]{}%
  }%
  \providecommand\includegraphics[2][]{%
    \GenericError{(gnuplot) \space\space\space\@spaces}{%
      Package graphicx or graphics not loaded%
    }{See the gnuplot documentation for explanation.%
    }{The gnuplot epslatex terminal needs graphicx.sty or graphics.sty.}%
    \renewcommand\includegraphics[2][]{}%
  }%
  \providecommand\rotatebox[2]{#2}%
  \@ifundefined{ifGPcolor}{%
    \newif\ifGPcolor
    \GPcolortrue
  }{}%
  \@ifundefined{ifGPblacktext}{%
    \newif\ifGPblacktext
    \GPblacktexttrue
  }{}%
  \let\gplgaddtomacro\g@addto@macro
  \gdef\gplbacktext{}%
  \gdef\gplfronttext{}%
  \makeatother
  \ifGPblacktext
    \def\colorrgb#1{}%
    \def\colorgray#1{}%
  \else
    \ifGPcolor
      \def\colorrgb#1{\color[rgb]{#1}}%
      \def\colorgray#1{\color[gray]{#1}}%
      \expandafter\def\csname LTw\endcsname{\color{white}}%
      \expandafter\def\csname LTb\endcsname{\color{black}}%
      \expandafter\def\csname LTa\endcsname{\color{black}}%
      \expandafter\def\csname LT0\endcsname{\color[rgb]{1,0,0}}%
      \expandafter\def\csname LT1\endcsname{\color[rgb]{0,1,0}}%
      \expandafter\def\csname LT2\endcsname{\color[rgb]{0,0,1}}%
      \expandafter\def\csname LT3\endcsname{\color[rgb]{1,0,1}}%
      \expandafter\def\csname LT4\endcsname{\color[rgb]{0,1,1}}%
      \expandafter\def\csname LT5\endcsname{\color[rgb]{1,1,0}}%
      \expandafter\def\csname LT6\endcsname{\color[rgb]{0,0,0}}%
      \expandafter\def\csname LT7\endcsname{\color[rgb]{1,0.3,0}}%
      \expandafter\def\csname LT8\endcsname{\color[rgb]{0.5,0.5,0.5}}%
    \else
      \def\colorrgb#1{\color{black}}%
      \def\colorgray#1{\color[gray]{#1}}%
      \expandafter\def\csname LTw\endcsname{\color{white}}%
      \expandafter\def\csname LTb\endcsname{\color{black}}%
      \expandafter\def\csname LTa\endcsname{\color{black}}%
      \expandafter\def\csname LT0\endcsname{\color{black}}%
      \expandafter\def\csname LT1\endcsname{\color{black}}%
      \expandafter\def\csname LT2\endcsname{\color{black}}%
      \expandafter\def\csname LT3\endcsname{\color{black}}%
      \expandafter\def\csname LT4\endcsname{\color{black}}%
      \expandafter\def\csname LT5\endcsname{\color{black}}%
      \expandafter\def\csname LT6\endcsname{\color{black}}%
      \expandafter\def\csname LT7\endcsname{\color{black}}%
      \expandafter\def\csname LT8\endcsname{\color{black}}%
    \fi
  \fi
    \setlength{\unitlength}{0.0500bp}%
    \ifx\gptboxheight\undefined%
      \newlength{\gptboxheight}%
      \newlength{\gptboxwidth}%
      \newsavebox{\gptboxtext}%
    \fi%
    \setlength{\fboxrule}{0.5pt}%
    \setlength{\fboxsep}{1pt}%
\begin{picture}(5040.00,3600.00)%
    \gplgaddtomacro\gplbacktext{%
      \csname LTb\endcsname
      \put(1078,660){\makebox(0,0)[r]{\strut{}1e-08}}%
      \put(1078,1000){\makebox(0,0)[r]{\strut{}1e-07}}%
      \put(1078,1340){\makebox(0,0)[r]{\strut{}1e-06}}%
      \put(1078,1680){\makebox(0,0)[r]{\strut{}1e-05}}%
      \put(1078,2020){\makebox(0,0)[r]{\strut{}1e-04}}%
      \put(1078,2359){\makebox(0,0)[r]{\strut{}1e-03}}%
      \put(1078,2699){\makebox(0,0)[r]{\strut{}1e-02}}%
      \put(1078,3039){\makebox(0,0)[r]{\strut{}1e-01}}%
      \put(1078,3379){\makebox(0,0)[r]{\strut{}1e+00}}%
      \put(1831,440){\makebox(0,0){\strut{}1e-04}}%
      \put(2591,440){\makebox(0,0){\strut{}1e-03}}%
      \put(3351,440){\makebox(0,0){\strut{}1e-02}}%
      \put(4112,440){\makebox(0,0){\strut{}1e-01}}%
      \put(3270,986){\makebox(0,0)[l]{\strut{}Slope 1}}%
    }%
    \gplgaddtomacro\gplfronttext{%
      \csname LTb\endcsname
      \put(341,2019){\rotatebox{-270}{\makebox(0,0){\strut{}$\|\sum q\|/\sum\|q\|$}}}%
      \put(2926,154){\makebox(0,0){\strut{}$\epsilon_{LS}$}}%
      \put(1870,3206){\makebox(0,0)[r]{\strut{}$\alpha=0$}}%
      \csname LTb\endcsname
      \put(1870,2986){\makebox(0,0)[r]{\strut{}$\alpha=5$}}%
      \csname LTb\endcsname
      \put(1870,2766){\makebox(0,0)[r]{\strut{}$\alpha=10$}}%
    }%
    \gplbacktext
    \put(0,0){\includegraphics[width={252.00bp},height={180.00bp}]{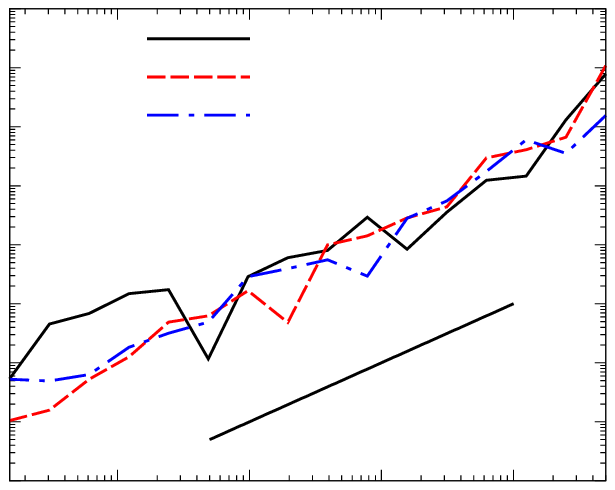}}%
    \gplfronttext
  \end{picture}%
\endgroup
   \caption{The imbalance between the total flow going in and out of the Glenn geometry (i.e., a measure of error in satisfying the conservation of mass) as a function of linear solver tolerance at three different Womersley numbers (solid; $\alpha=0$, dashed; $\alpha=5$, and dashed-dot; $\alpha=10$). }
  \label{fig:ser}    
\end{figure}

Since there is no closed-form solution available for the assessment of the accuracy of this case, the net flow through three branches (SVC, LPA, and RPA) is computed to obtain a measure of the error in the solution of the continuity equation (Figure \ref{fig:ser}). 
It is observed that this error is linearly proportional to the tolerance by which the linear system of equations is solved.   
Varying the Womersley number has little effect on this error, demonstrating the robustness of the proposed method for a wide range of oscillatory flows. 

In terms of cost, it takes a full day and 3,000 CPU-hours to simulate this case on a 128 processor cluster using a conventional CFD solver \cite{meng2020time}. 
These numbers are reduced to 5.75 seconds and 0.2 CPU-hours when using the proposed formulation (considering the most expensive case of $\alpha =5$). 
This four orders of magnitude reduction in wall clock time and computational cost is roughly equal to the number of time steps for the conventional solver, which is 10,000 (to simulate a total of 5 cycles to ensure cycle-to-cycle convergence).
As stated earlier, this large gap in wall-clock time is independent of the number of computed modes as these computations are embarrassingly parallelizable across different modes. 

\section{Conclusions} \label{sec:conclusion}
A new technique was proposed for solving the Stokes equations in the frequency domain. 
This method only relied on real arithmetics and allowed for the use of similar order shape functions for pressure and velocity unknowns. 
The method was built around the addition of a stabilization term in the form of the Laplacian of pressure to the continuity equation (Eq. \eqref{stokes_S2}). 
This stabilization term that was derived systematically from the momentum equation produced a complex-valued stabilization parameter $\tau$ (Eq. \eqref{tau}).
Varying the single adjustable parameter that appeared in its definition showed that it had little effect on the accuracy of the overall scheme. 
The results also showed the robustness of the proposed method against variation in the Womersley number, mesh size, and geometrical complexity. 
The proposed method was very computationally efficient, enabling a typical simulation, which otherwise takes hours to simulate using a conventional solver, to be simulated in seconds.


\begin{thebibliography}{10}

\bibitem{meng2020time}
Chenwei Meng, Anirban Bhattacharjee, and Mahdi Esmaily.
\newblock A scalable spectral {S}tokes solver for simulation of time-periodic
  flows in complex geometries.
\newblock {\em Journal of Computational Physics}, 445:110601, 2021.

\bibitem{ladyzhenskaya1969mathematical}
Olga~Alexandrovna Ladyzhenskaya.
\newblock {\em The mathematical theory of viscous incompressible flow},
  volume~2.
\newblock Gordon and Breach New York, 1969.

\bibitem{babuvska1971error}
Ivo Babuska.
\newblock Error-bounds for finite element method.
\newblock {\em Numerische Mathematik}, 16(4):322--333, 1971.

\bibitem{brezzi1974existence}
Franco Brezzi.
\newblock On the existence, uniqueness and approximation of saddle-point
  problems arising from {L}agrangian multipliers.
\newblock {\em ESAIM: Mathematical Modelling and Numerical
  Analysis-Mod{\'e}lisation Math{\'e}matique et Analyse Num{\'e}rique},
  8(R2):129--151, 1974.

\bibitem{bazilevs2013computational}
Yuri Bazilevs, Kenji Takizawa, and Tayfun~E Tezduyar.
\newblock {\em Computational fluid-structure interaction: methods and
  applications}.
\newblock John Wiley \& Sons, 2013.

\bibitem{Esmaily2015DS}
Mahdi Esmaily, Yuri Bazilevs, and Alison Marsden.
\newblock Impact of data distribution on the parallel performance of iterative
  linear solvers with emphasis on {CFD} of incompressible flows.
\newblock {\em Computational Mechanics}, 55(1):93--103, 2015.

\bibitem{saad1986gmres}
Youcef Saad and Martin~H Schultz.
\newblock {GMRES}: A generalized minimal residual algorithm for solving
  nonsymmetric linear systems.
\newblock {\em SIAM Journal on scientific and statistical computing},
  7(3):856--869, 1986.

\bibitem{Esmaily2015BIPN}
Mahdi Esmaily, Yuri Bazilevs, and Alison Marsden.
\newblock A bi-partitioned iterative algorithm for solving linear systems
  arising from incompressible flow problems.
\newblock {\em Computer Methods in Applied Mechanics and Engineering},
  286(1):40--62, 2015.

\bibitem{shakib1991new}
Farzin Shakib, Thomas~JR Hughes, and Zden{\v{e}}k Johan.
\newblock A new finite element formulation for computational fluid dynamics: X.
  the compressible euler and navier-stokes equations.
\newblock {\em Computer Methods in Applied Mechanics and Engineering},
  89(1-3):141--219, 1991.

\bibitem{bazilevs2007variational}
Y~Bazilevs, VM~Calo, JA~Cottrell, TJR Hughes, A~Reali, and G~Scovazzi.
\newblock Variational multiscale residual-based turbulence modeling for large
  eddy simulation of incompressible flows.
\newblock {\em Computer Methods in Applied Mechanics and Engineering},
  197(1-4):173--201, 2007.

\bibitem{METIS}
George Karypis and Vipin Kumar.
\newblock {MeTis: Unstructured Graph Partitioning and Sparse Matrix Ordering
  System, Version 4.0}.
\newblock \url{http://www.cs.umn.edu/~metis}, 2009.

\bibitem{Womersley1955method}
John~R Womersley.
\newblock Method for the calculation of velocity, rate of flow and viscous drag
  in arteries when the pressure gradient is known.
\newblock {\em The Journal of Physiology}, 127(3):553, 1955.

\bibitem{arbia2014numerical}
Gregory Arbia, Chiara Corsini, Mahdi Esmaily, Alison~L Marsden, Francesco
  Migliavacca, Giancarlo Pennati, Tain-Yen Hsia, Irene~E Vignon-Clementel,
  Modeling of~Congenital Hearts Alliance (MOCHA)~Investigators, et~al.
\newblock Numerical blood flow simulation in surgical corrections: what do we
  need for an accurate analysis{?}
\newblock {\em Journal of Surgical Research}, 186(1):44--55, 2014.

\end{thebibliography}
\end{document}